\documentclass[reqno]{amsart}
\usepackage{amssymb}
\usepackage[usenames, dvipsnames]{color}
\usepackage{txfonts}
\usepackage{verbatim}
\usepackage{marginnote}
\usepackage{todonotes}
\usepackage[normalem]{ulem}
\usepackage{cancel}
\usepackage{soul}

\theoremstyle{plain}
\newtheorem{theorem}{Theorem}[section]
\newtheorem{lemma}[theorem]{Lemma}

\newtheorem{proposition}[theorem]{Proposition}

\theoremstyle{definition}

\theoremstyle{remark}
\newtheorem{remark}[theorem]{Remark}

\newcommand{\bR}{{\mathbb R}}
\newcommand{\bN}{{\mathbb N}}

\def\La{\Lambda}
\def\t{\tilde}
\def\th{\theta}

\def\Dl{\Delta}
\def\lt{\left}
\def\rt{\right}

\def\i{\infty}

\def \ls{\lesssim}
\def\p{\partial}
\def\f{\frac}
\def\na{\nabla}

\def\O{\Omega}

\def\q{\quad}
\def\qq{\qquad}

\def\be{\begin{equation}}
\def\ee{\end{equation}}
\def\bes{\begin{equation*}}
\def\ees{\end{equation*}}
\def\bali{\begin{aligned}}
\def\eali{\end{aligned}}
\def\ML{\mathcal{L}}
\def\pa{\|}

\numberwithin{equation}{section}

\makeatletter
\def\dashint{\operatorname%
{\,\,\text{\bf--}\kern-.98em\DOTSI\intop\ilimits@\!\!}}
\makeatother

\begin{document}


\title[Global regularity of MHD-Boussinesq]{Global regularity of solutions for the 3D non-resistive and non-diffusive MHD-Boussinesq system with axisymmetric data}

\author[X. Pan]{Xinghong Pan}
\address[X. Pan]{Department of Mathematics, Nanjing University of Aeronautics and Astronautics, Nanjing 211106, China}

\email{xinghong\_87@nuaa.edu.cn}

\thanks{X. Pan is supported by Natural Science Foundation of Jiangsu Province (No. SBK2018041027), Double Innovation Scheme of Jiangsu Province and National Natural Science Foundation of China (No. 11801268).}



\begin{abstract}
In this paper, we will show that solutions of the three-dimensional non-resistive and non-diffusive MHD-Boussinesq system are globally regular if the initial data is axisymmetric and the swirl components of the velocity and the magnetic vorticity are zero. Our main result extends previous ones on the three-dimensional non-resistive MHD system and non-diffusive Boussinesq system, and the method used here can also be applied to the magnetic Rayleigh-B\'{e}nard convection system.


\end{abstract}
\maketitle

\section{Introduction}
In this paper, we consider the global regularity problem for the three-dimensional (3D) magnetohydrodynamics (MHD) -Boussinesq system
\begin{equation}\label{MB}
\left\{
\begin{aligned}
&
\p_tu+u\cdot\nabla u+\nabla p-\mu\Delta u=h\cdot\nabla h+\rho e_3,\\
&\p_t h+u\cdot\nabla h-h\cdot\nabla u-\nu\Dl h=0,\\
&\p_t\rho+u\cdot\na\rho-\kappa\Dl\rho=0,\\
&\nabla\cdot u=\nabla\cdot h=0.\\
\end{aligned}
\right.
\end{equation}
Here $u(t,x),\ h(t,x)\in\bR^3,p(t,x)\in\bR$ and $\rho(t,x)\in\bR$ represent the velocity, magnetic field, pressure and temperature fluctuation. The vector $e_3=(0,0,1)$ is the unit vector in the vertical direction. $\mu\geq 0,\nu\geq 0$ and $\kappa\geq 0$ stand for the constant viscosity, magnetic resistivity  and thermal diffusivity, respectively. The MHD-Boussinesq system models the convection of an incompressible flow driven by the buoyant effect of a thermal field and the Lorenz force, generated by the magnetic field.


We say that the MHD-Boussinesq system is non-resistive and non-diffusive, which means $\mu>0$, but $\nu=\kappa=0$. Without loss of generality, we set $\mu=1$ and system \eqref{MB} becomes
\begin{equation}\label{MB1}
\left\{
\begin{aligned}
&
\p_tu+u\cdot\nabla u+\nabla p-\Delta u=h\cdot\nabla h+\rho e_3,\\
&\p_t h+u\cdot\nabla h-h\cdot\nabla u=0,\\
&\p_t\rho+u\cdot\na\rho=0,\\
&\nabla\cdot u=\nabla\cdot h=0.\\
\end{aligned}
\right.
\end{equation}

The local well-posedness result of \eqref{MB1} can be founded in \cite{LP:JDE2017}. However, the global well-posedness is still wildly open even for the Navier-Stokes equations ($h=\rho\equiv 0$), let alone for the system \eqref{MB1}. In this paper, we will show that a family of axisymmetic solutions to \eqref{MB1} are globally as regular as their initial data.

In the following, we will carry out our proof in the cylindrical coordinates $(r,\,\th,\,z)$. That is, for $x=(x_1,\,x_2,\,x_3)\in\mathbb{R}^3$
\bes
r=\sqrt{x_1^2+x_2^2},\q \th=\arctan\frac{x_2}{x_1},\q z=x_3.
\ees
And the axisymmetic solution of system \eqref{MB1} is given by
\[
u=u^r(t,r,z)e_r+u^{\th}(t,r,z)e_{\th}+u^z(t,r,z)e_z,
\]
\[
h=h^r(t,r,z)e_r+h^{\th}(t,r,z)e_{\th}+h^z(t,r,z)e_z,
\]
\[
\rho=\rho(t,r,z),
\]
where the basis vectors $e_r,e_\th,e_z$ are
\[
e_r=(\frac{x_1}{r},\frac{x_2}{r},0),\quad e_\th=(-\frac{x_2}{r},\frac{x_1}{r},0),\quad e_z=(0,0,1).
\]
We will prove the global regularity of the following family of axisymmetric solutions
\be\label{assolution}
u=u^r(t,r,z)e_r+u^z(t,r,z)e_z,\q h=h^\th(t,r,z) e_\th,\q \rho=\rho(t,r,z).
\ee
Denote
\bes
\Phi_{k,c}(t):=c\underbrace{\exp(\cdots\exp}_{k~\text{times}}(ct)\cdots).
\ees More precisely, we have the following theorem.

\begin{theorem}\label{th1}
Let $u_0$, $h_0$ and $\rho_0$ be all axially symmetric data with $\na\cdot u_0=0$. Besides, we assume that $u^\th_0=h^r_0=h^z_0=0$. If $(u_0,h_0,\rho_0)\in H^2(\bR^3)$ and $H_0:=\f{h^\th_0}{r}\in L^\i(\bR^3)$, then there exists a unique global solution $(u,h,\rho)$ to the MHD-Boussinesq system \eqref{MB1} with data $(u_0,h_0,\rho_0)$, which satisfies
\be\label{energy}
\|(u,h,\rho)(t,\cdot)\|^2_{H^2}+\int^t_0\|\na u(t,\cdot)\|^2_{H^2}ds\leq \Phi_{3,c_0}(t),
\ee
where $c_0$ is a positive constant depending only on $H^2$ norms of $u_0, h_0, \rho_0$ and $L^\i$ norm of $ H_0$.
\end{theorem}
\begin{remark}
It is not hard to extend the result of Theorem \eqref{th1} to the case where $\mu>0$, $\nu\geq 0$ and $\kappa\geq 0$ in \eqref{MB} with the same initial data as that in Theorem \eqref{th1}.
\end{remark}
\qed
\begin{remark}
When $h^\th\equiv 0$, the global well-posedness result for the axisymmetric Navier-Stokes-Boussinesq can be found in \cite{AHK:2011DCDS,HR:2010AIHP}. While if $\rho\equiv0$, see \cite{Lz:2015JDE} for the global well-posedness result for the axisymmetric MHD system. Our main result can be viewed as an extension of those in the above papers.
\end{remark}
\qed
\begin{remark}
 Define
\bes
 H:=\f{h^\th}{r},\q \O:=\f{w^\th}{r},\q w^\th=\p_z u^r-\p_r u^z.
\ees

 The proof of Theorem \ref{th1} strongly depends on the special structure of the MHD-Boussinesq system in axisymmtric case with zero swirl components of the velocity and the magnetic vorticity. We will show that $ H$ and $\rho$ satisfy the same transport equations and $\O$ satisfies a linear diffusive equation with inhomogeneous terms involving only in $ H$ and $\rho$. See \eqref{MB2}. Then the $L^\i_tL^2_x$ norm of $\O$ will be obtained. This is a key step for us to bootstrap the regularity of $u,\ h$ and $\rho$.

 Our proof combines the ideas that in \cite{HR:2010AIHP} and \cite{Lz:2015JDE}. Here we outline the main differences. Compared with that in \cite{HR:2010AIHP},  we need to deal with the extra term $\p_z H$ in \eqref{MB2} and later much more estimates on the magnetic filed $h^\th e_\th$ are needed, which are nontrivial. Compared with that in \cite{Lz:2015JDE}, in our paper, the $L^\i_t L^2_x\cap L^2_t H^1_x$ of $\O$ can not be obtained from the system \eqref{MB2} due to the appearance of $\f{\p_r\rho}{r}$. So the estimate $\|u^r/r\|_{L^1_t L^\i_x}$ in \cite[Lemma 2.2]{Lz:2015JDE} is not applicable to us.
\end{remark}
\qed
\begin{remark}
This MHD-Boussinesq system \eqref{MB} is closely related to a type of the Rayleigh-B\'{e}nard convection, which occurs in a horizontal layer of conductive fluid heated from below, with a presence of a magnetic field. The only difference between the magnetic Rayleigh-B\'{e}nard convection system and the MHD-Boussinesq system is that $\eqref{MB}_3$ is replaced by the following equation
\bes
\p_t\rho+u\cdot\na\rho-\kappa\Dl\rho=u^3.
\ees
Various physical theories and numerical experiments have been developed to study the magnetic Rayleigh-B\'{e}nard convection and related equations. See, for example, \cite{MR:ARMA2003,TM:PD1981} and references therein.  The result in Theorem \ref{th1} can also be applied to the following non-resistive and non-diffusive magnetic Rayleigh-B\'{e}nard convection system
\bes
\left\{
\begin{aligned}
&
\p_tu+u\cdot\nabla u+\nabla p-\Delta u=h\cdot\nabla h+\rho e_3,\\
&\p_t h+u\cdot\nabla h-h\cdot\nabla u=0,\\
&\p_t\rho+u\cdot\na\rho=u^3,\\
&\nabla\cdot u=\nabla\cdot h=0.\\
\end{aligned}
\right.
\ees
The proof is essentially the same as that for \eqref{MB1} with little difference. We omit the details.
\end{remark}
\qed

If the fluid is not affected by the temperature, then our system \eqref{MB} is reduced to the classical MHD system. There already have been many studies and fruitful results related to the well-posedness of the MHD system. Sermange-Temam \cite{ST:1983CPAM} established the local existence and uniqueness of the solution and particularly the 2D local strong solution was proved to be global. Cao et al. in \cite{CW:2011ADVANCE,CRW:2013JDE} proved the global regularity of the MHD system for a variety of combinations of partial dissipation and diffusion in 2D space. Lin-Xu-Zhang \cite{LXZ:JDE2015} proved the global well-posedness of classical solutions for 2D non-resistive MHD under the assumption that the initial data is a small perturbation of a nonzero constant magnetic field. See also \cite{RWXZ:2014JFA} for similar results. For the 3D case, readers can see \cite{LZ:CPAM2014,XZ:2015JMA} for related results. Cai-Lei \cite{CL:2018ARMA} and He-Xu-Yu \cite{HXY:2018ANNPDE} proved the global well-posedness of small initial data for the idea (inviscid and non-resistive) MHD system. Lei \cite{Lz:2015JDE} proved the global regularity of classical solutions to the 3D MHD system with a family of axisymmetric large data. We also emphasized some partial regularity results and blow up criteria in \cite{HX:2005JFA,HX:2005JDE,CMZ:2008CMP,LZ:2009DCDS} and references therein.

 On the other hand, if the fluid is not affected by the Lorentz force, then our system \eqref{MB1} is the classical Boussinesq system without diffusion. Many works and efforts have been made to study the well-posedness of the Cauchy problem for the Boussinesq system. In 2D case, Chae \cite{Chae:2006ADV} and Hou-Li \cite{HL:2005DCDS} independently proved the global regularity of solutions to the 2D Boussinesq system. And also Chae \cite{Chae:2006ADV} considered the case of zero viscosity and non-zero diffusion. See \cite{AH:2007JDE,HK:2009INDIANA} for related results in critical space. For 3D case, Abidi et al. \cite{AHK:2011DCDS} and Hmidi-Rousset \cite{HR:2010AIHP,HR:2011JFA} proved the global well-posedness of the Cauchy problem for the 3D axisymmetric Boussinesq system without swirl. Readers can see \cite{LLT:2013JDE,CW:2013ARMA} and references therein for more regularity results on the Boussinesq system.

 For the full MHD-Boussinesq system, there are also some works concentrated on the global well-posedness of weak and strong solutions. See \cite{BG:2016JDE,BL:2017JDE} and references therein for 2D cases. In the 3D case, Larios-Pei \cite{LP:JDE2017} proved the local well-posedness results in Sobolev space. Liu-Bian-Pu \cite{LBP:2019ZAMP} proved the global well-posedness of strong solutions with nonlinear damping term in the momentum equations. Recently, Bian-Pu \cite{BP:2020JMFM} proved the global regularity of a family of axially symmetric large solutions to the MHDB system without magnetic resistivity and thermal diffusivity under the assumption that the support of the initial thermal fluctuation is away from the $z$-axis and its projection on to the $z$-axis is compact. In this paper, we will improve the result in \cite{BP:2020JMFM} by removing the ``support set" assumption on the data of the thermal fluctuation. Regarding the MHD-B\'{e}nard system, some progress has also been made in 2D and 3D cases. See, e.g., \cite{ZFN:2013AML,CD:2015JMFM,Yk:2017MMAS,ZT:2018AA} and references therein.

Our paper is organized as follows. In Section 2, we reformulate our system in cylindrical coordinates and prove an a priori $L^\i_tL^2_x$ estimate for $\O$. In Section 3, we give the $H^1$ a priori estimate of the solution. In Section 4, we give the $H^2$ a priori estimate of the solution and prove Theorem \ref{th1}. Throughout the paper, we use $C$ or $c$ to denote a generic constant which may be different from line to line. We also apply $A\lesssim B$ to denote $A\leq CB$.

\section{Reformulation of the system and ${L^\i_tL^2_x}$ estimate of ${\O}$}

The axisymmetric MHD-Boussinesq system \eqref{MB1} in cylindrical coordinates read
\begin{equation}\label{MBAS}
\left\{
\begin{aligned}
&\p_tu^r+(u^r\p_r+u^z\p_z)u^r -\frac{(u^\th)^2}{r}+\p_r P=(h^r\p_r+h^z\p_z)h^r-\frac{(h^\th)^2}{r}+(\Delta-\frac{1}{r^2})u^r, \\
&\p_tu^\th+(u^r\p_r+u^z\p_z) u^\th+\frac{u^\th u^r}{r}=(h^r\p_r+h^z\p_z)h^\th+\frac{h^rh^\th}{r}+(\Delta-\frac{1}{r^2})u^\th , \\
&\p_tu^z+(u^r\p_r+u^z\p_z)u^z+\p_z P=(h^r\p_r+h^z\p_z)h^z+\Delta u^z+\rho ,                                    \\
&\p_th^r+(u^r\p_r+u^z\p_z)h^r-(h^r\p_r+h^z\p_z)u^r=0,\\
&\p_th^\th+(u^r\p_r+u^z\p_z)h^\th-(h^r\p_r+h^z\p_z)u^\th+\frac{u^\th h^r}{r}-\frac{h^\th u^r}{r}=0,\\
&\p_th^z+(u^r\p_r+u^z\p_z)h^z-(h^r\p_r+h^z\p_z)u^z=0,\\
&\p_t \rho+(u^r\p_r+u^z\p_z)\rho=0,\\
&\nabla\cdot u=\p_ru^r+\frac{u^r}{r}+\p_zu^z=0,\quad \nabla\cdot h=\p_rh^r+\frac{h^r}{r}+\p_zh^z=0,
\end{aligned}
\right.
\end{equation}
where the pressure $P=p+\f{1}{2}|h|^2$ and $\Delta=\frac{\p^2}{\p r^2}+\frac{1}{r}\frac{\p}{\p r}+\frac{\p^2}{\p z^2}$ is the usual Laplacian operator. By the uniqueness of local solutions, it is easy to see that if the initial data satisfy $u^\th_0=h^r_0=h^z_0=0$, then the solution of \eqref{MBAS} will be the form of \eqref{assolution}. In this situation, \eqref{MBAS} can be simplified as
\begin{equation}\label{MB4}
\left\{
\begin{aligned}
&\p_t u^r+(u^r\p_r+u^z\p_z)u^r+\p_r P=-\frac{(h^\th)^2}{r}+(\Delta-\frac{1}{r^2})u^r, \\
&\p_t u^z+(u^r\p_r+u^z\p_z)u^z+\p_z P=\Delta u^z+\rho ,                                    \\
&\p_t h^\th+(u^r\p_r+u^z\p_z)h^\th-\frac{u^r}{r}h^\th=0,\\
&\p_t \rho+(u^r\p_r+u^z\p_z)\rho=0,\\
&\f{1}{r}\p_r(ru^r)+\p_zu^z=0.
\end{aligned}
\right.
\end{equation}
Denote $ H:=\f{h^\th}{r}$ and $\O:=\f{w^\th}{r}$. From \eqref{MB4}, we can get
\be\label{MB2}
\lt\{
\bali
&\p_t\O+u\cdot\na \O=(\Dl+\f{2}{r}\p_r)\O-\p_z H^2-\f{\p_r\rho}{r},\\
&\p_t H+u\cdot\na H=0,\\
&\p_t \rho+u\cdot\na \rho=0.
\eali
\rt.
\ee
First we have the following Proposition.
\begin{proposition}\label{pl2}
Let $(u,h,\rho)$ be a smooth solution of \eqref{MB4}, then we have\\
(1) for $p\in[1,\i]$ and $t\in\bR_+$, we have
\be\label{jrho}
\|( H(t),\rho(t))\|_{L^p}\leq \|( H_0,\rho_0)\|_{L^p};
\ee
(2) for $u_0,h_0, \rho_0\in L^2$ and $t\in\bR_+$, we have
\be\label{l2estimate}
\|(u(t),h(t))\|^2_{L^2}+\int^t_0\|\na u(s)\|ds\leq C_0(1+t)^2,
\ee
where $C_0$ depends only on $\|(u_0,h_0)\|_{L^2}$ and $\|\rho_0\|_{L^2}$.
\end{proposition}

\noindent\textbf{Proof of Proposition \ref{pl2}.}
\begin{proof}
 The estimate in \eqref{jrho} is classical for the transport equation with finite $p$. While if $p=\i$, it is just the maximum principle. For the estimate in \eqref{l2estimate}, we proceed the standard $L^2$ inner product estimate of system \eqref{MB1}. Then we have
 \be\label{l21}
 \f{1}{2}\f{d}{dt}\|(u(t),h(t))\|^2_{L^2}+\|\na u(t)\|^2_{L^2}\leq \|u(t)\|_{L^2}\|\rho(t)\|_{L^2}.
 \ee
This indicates that
\bes
\f{d}{dt}\|(u(t),h(t))\|_{L^2}\leq 2\|\rho(t)\|_{L^2}.
\ees
Integration on time indicates that
\bes
\bali
\|(u(t),h(t))\|_{L^2}&\leq \|(u_0,h_0)\|_{L^2}+2\int^t_0\|\rho(\tau)\|_{L^2}d\tau\\
                     &\leq \|(u_0,h_0)\|_{L^2}+2\|\rho_0\|_{L^2}t.
\eali
\ees
Inserting this into \eqref{l21} and integration on time, we have
\bes
\bali
&\q\f{1}{2}\|(u(t),h(t))\|^2_{L^2}+\int^t_0\|\na u(s)\|^2_{L^2}ds\\
&\leq \f{1}{2}\|(u_0,h_0)\|^2_{L^2}+\big(\|(u_0,h_0)\|_{L^2}+2\|\rho_0\|_{L^2}t\big)\|\rho_0\|_{L^2} t.
\eali
\ees
This gives \eqref{l2estimate}.
\end{proof}

Based on Proposition \ref{pl2}, we have the following Proposition which gives the a priori $L^\i_tL^2_x$ estimate of $\O$.
\begin{proposition}\label{pomega}
Suppose $(u,h,\rho)$ be the smooth solution of \eqref{MB1} with initial data $(u_0,h_0,\rho_0)$ satisfying assumptions in Theorem \ref{th1}, then we have, for $t\in\bR_+$,
\be\label{eomega}
\|\O(t)\|_{L^2}\leq \Phi_{1,c_0}(t),
\ee
where $c_0$ is a positive constant depending only on $H^2$ norms of $u_0, h_0, \rho_0$ and $L^\i$ norm of $ H_0$.
\end{proposition}
\noindent Before proving Proposition \ref{pomega}, we collect some useful estimates and identities.
\begin{lemma}[Proposition 3.1, 3.2 and Lemma 3.3 of \cite{HR:2010AIHP}]
Denote $\ML=(\Dl+\f{2}{r}\p_r)^{-1}\f{\p_r}{r}$ and $\t{\ML}=(\Dl+\f{2}{r}\p_r)^{-1}\f{\p_z}{r}$. Suppose $\rho\in H^2(\bR^3)$ be axisymmetric, then for every $p\in[2,+\i)$, there exists an absolute constant $C_p>0$ such that
\be\label{OL}
\|\ML\rho\|_{L^p}\leq C_p\|\rho\|_{L^p},\q \|\t{\ML}\rho\|_{L^p}\leq C_p\|\rho\|_{L^p}.
\ee
Moreover, for any smooth axisymmetric function $f$, we have the identity
\be\label{OL1}
\ML\p_rf=\f{f}{r}-\ML\Big(\f{f}{r}\Big)-\p_z\t{\ML}f.
\ee
\end{lemma}

\begin{lemma}
For $1<p<+\i$, there exists an absolute constant $C_p>0$ such that
\be\label{Omega}
\|\na\f{u^r}{r}\|_{L^p}\leq C_p\|\O\|_{L^p}.
\ee
\end{lemma}
The proof of this lemma can be founded in many literatures , such as \cite[A.5 on page 3213]{Lz:2015JDE}, \cite[Lemma 2.3]{CFZ:2017DCDS} or \cite[Proposition 2.5]{MZ:2013CMP}.

\noindent\textbf{Proof of Proposition \ref{pomega}}
\begin{proof}
Applying $\ML$ to $\eqref{MB2}_3$, we get
\be\label{rho}
\p_t\ML\rho+u\cdot\na\ML\rho=-[\ML,u\cdot\na]\rho,
\ee
where $[A,B]=AB-BA$ is the commutator.

Denote $L:=\O-\ML\rho$. Subtracting \eqref{rho} from $\eqref{MB2}_1$, we have
\be\label{GAMMA}
\p_t L+u\cdot\na L-(\Dl+\f{2}{r}\p_r) L=[\ML,u\cdot\na]\rho-\p_z  H^2.
\ee
Taking $L^2$ inner product of \eqref{GAMMA}, using integration by parts and divergence-free condition of $u$, we get
\bes
\bali
&\q\f{1}{2}\f{d}{dt}\pa L(t)\pa^2_{L^2}+\pa\na L(t)\pa^2_{L^2}\\
&\leq\int_{\bR^3}\ML(u\cdot\na\rho)L dx-\int_{\bR^3}u\cdot\na(\ML\rho) L dx-\int_{\bR^3}\p_z H^2 L dx\\
&\leq\int_{\bR^3}\ML(u\cdot\na\rho) L dx+\int_{\bR^3}(\ML\rho) u\cdot\na L dx+\int_{\bR^3} H^2\p_z L dx\\
&:=I_1+I_2+I_3.
\eali
\ees
Next we will estimate $I_i\ (i=1,2,3)$ term by term. For $I_1$, first we make some computation on $\ML(u\cdot\na\rho)$.
\bes
\bali
\ML(u\cdot\na\rho)&=\ML(\na\cdot(u\rho))\\
                  &=\ML\Big(\p_r(u^r\rho)+\f{1}{r}(u^r\rho)+\p_z(u^z\rho)\Big).
\eali
\ees
From \eqref{OL1}, we have
\bes
\bali
\ML(u\cdot\na\rho)=&\ML\p_r(u^r\rho)+\ML\Big(\f{u^r\rho}{r}\Big)+\ML\p_z(u^z\rho)\\
=&\f{u^r}{r}\rho-\p_z\t{\ML}(u^r\rho)+\p_z\ML(u^z\rho),
\eali
\ees
where we have used the fact that $\p_z$ is commutated with $\ML$.
\\
Then, using integration by parts, we get
\bes
\bali
I_1&=\int_{\bR^3}\f{u^r}{r}\rho L dx+\int_{\bR^3}\t{\ML}(u^r\rho)\p_z L dx-\int_{\bR^3}\ML(u^z\rho)\p_z Ldx\\
   &=I^1_1+I^2_1+I^3_1.
\eali
\ees
Using H\"{o}lder inequality, Sobolev embedding and \eqref{Omega}, we have
\bes
\bali
|I^1_1|&\leq\pa\f{u^r}{r}\pa_{L^6}\pa\rho\pa_{L^3}\pa L\pa_{L^2}\\
       &\leq \pa\na\f{u^r}{r}\pa_{L^2}\pa\rho\pa_{L^3}\pa L\pa_{L^2}\\
       &\leq \pa\O\pa_{L^2}\pa\rho\pa_{L^3}\pa L\pa_{L^2}\\
       &\leq(\pa L\pa_{L^2}+\pa\ML\rho\pa_{L^2})\pa\rho\pa_{L^3}\pa L\pa_{L^2}.
\eali
\ees
Using \eqref{OL}, \eqref{jrho} and Sobolev embedding, we have
\bes
\bali
|I^1_1| &\leq C(\pa L\pa_{L^2}+\pa\rho\pa_{L^2})\pa\rho\pa_{L^3}\pa L\pa_{L^2}\\
        &\leq C\pa\rho_0\pa_{L^3}\pa L\pa^2_{L^2}+C\pa\rho_0\pa_{L^2}\pa\rho_0\pa_{L^3}\pa L\pa_{L^2}\\
       &\leq C\pa\rho_0\pa_{H^2}\pa L\pa^2_{L^2}+C\pa\rho_0\pa^2_{H^2}\pa L\pa_{L^2}\\
       &\leq C\big(\pa\rho_0\pa_{H^2}+1\big)\pa L\pa^2_{L^2}+C\pa\rho_0\pa^4_{H^2}.
\eali
\ees
From \eqref{OL}, Proposition \ref{pl2} and using H\"{o}lder inequality, Young inequality, we have
\bes
\bali
&\q|I^2_1|+|I^3_1|\\
&\leq\Big(\pa\t{\ML}(u^r\rho)\pa_{L^2}+\pa{\ML}(u^r\rho)\pa_{L^2}\Big)\pa\p_z L\pa_{L^2}\\
       &\leq C\pa u^r\rho\pa_{L^2}\pa\p_z L\pa_{L^2}\\
       &\leq C\pa\rho_0\pa_{L^\i}\pa u\pa_{L^2}\pa\p_z L\pa_{L^2}\\
       &\leq  C\pa\rho_0\pa^2_{L^\i}\pa u\pa^2_{L^2}+\f{1}{4}\pa\p_z L\pa^2_{L^2}\\
       &\leq  C_0(1+t)^2+\f{1}{4}\pa\p_z L\pa^2_{L^2},
\eali
\ees
where $C_0$ is a positive constant depending only on $H^2$ norms of $u_0, h_0, \rho_0$ and $L^\i$ norm of $ H_0$. Also, the same techniques as above imply
\bes
\bali
&\q|I^2|+|I^3|\\
&\leq\Big(\pa(\ML\rho) u\pa_{L^2}+\pa  H^2\pa_{L^2}\Big)\pa\na L\pa_{L^2}\\
&\leq\Big(\pa\ML\rho\pa_{L^3}\pa u\pa_{L^6}+\pa  H\pa_{L^\i}\pa H\pa_{L^2}\Big)\pa\na L\pa_{L^2}\\
&\leq\Big(\pa\rho\pa_{L^3}\pa \na u\pa_{L^2}+\pa  H_0\pa_{L^\i}\pa H_0\pa_{L^2}\Big)\pa\na L\pa_{L^2}\\
&\leq \Big(\pa\rho_0\pa_{L^3}\pa \na u\pa_{L^2}+\pa H_0\pa_{L^\i}\pa h_0\pa_{H^2}\Big)^2+\f{1}{4}\pa\na L\pa^2_{L^2}\\
&\leq C_0\Big(1+\pa\na u\pa^2_{L^2}\Big)+\f{1}{4}\pa\na L\pa^2_{L^2}.
\eali
\ees
The above estimates indicate that
\bes
\bali
&\q\f{d}{dt}\pa L(t)\pa^2_{L^2}+\pa\na L(t)\pa^2_{L^2}\\
&\leq C_0\Big(1+\pa\na u\pa^2_{L^2}\Big)+C_0(1+t)^2\\
&\q+C\big(\pa\rho_0\pa_{H^2}+1\big)\pa L\pa^2_{L^2}+C\pa\rho_0\pa^4_{H^2}.
\eali
\ees
Gronwall inequality indicates that
\bes
\pa L(t)\pa^2_{L^2}+\int^t_0\pa\na L(s)\pa^2_{L^2}ds\leq \Phi_{1,c_0}(t).
\ees
Then we have
\bes
\bali
\pa\O(t)\pa_{L^2}&\leq \pa L\pa_{L^2}+\pa\ML\rho\pa_{L^2}\\
&\leq \pa L\pa_{L^2}+C\pa\rho\pa_{L^2}\\
&\leq \pa L\pa_{L^2}+\pa\rho_0\pa_{L^2}\leq\Phi_{1,c_0}(t).
\eali
\ees
This proves Proposition \ref{pomega} and \eqref{eomega} is valid.
\end{proof}
\section {${H^1}$ estimate of the solution}
In this section, we give a prior $H^1$ estimate for the solution of system \ref{MB4}. We have the following Proposition.
\begin{proposition}\label{pH1}
Suppose $(u,h,\rho)$ be the smooth solution of \eqref{MB1} with initial data $(u_0,h_0,\rho_0)$ satisfying assumptions in Theorem \ref{th1}, then we have, for $t\in\bR_+$,
\be\label{h1}
\bali
&\q\pa (\na u(t),\na h(t),\na\rho(t)) \pa^2_{L^2}+\int^t_0\pa \na^2 u(s) \pa^2_{ L^2}ds\leq \Phi_{2,c_0}(t),
\eali
\ee
where $c_0$ is a positive constant depending only on $H^2$ norms of $u_0, h_0, \rho_0$ and $L^\i$ norm of $ H_0$.
\end{proposition}

\subsection{${L^\i_tL^2}\cap {L^2_tH^1} $ estimate of $\na u$}

In cylindrical coordinates, the vorticity of the swirl-free axisymmetric velocity $u$ is given by $w=\na\times u=w^\th e_\th$ and $w^\th$ satisfies
\bes
\p_t w^\th+u\cdot\na w^\th-(\Dl-\f{1}{r^2})w^\th-\f{u^r}{r}w^\th=-\p_z\f{(h^\th)^2}{r}-\p_r\rho.
\ees
Performing the standard $L^2$ inner product, we have
\bes
\bali
&\q\f{1}{2}\f{d}{dt}\pa w^\th \pa^2_{L^2}+\pa \na w^\th \pa^2_{L^2}+\big\pa \f{w^\th}{r} \big\pa^2_{L^2}\\
&\leq \int_{\bR^3}\f{u^r}{r}(w^\th)^2dx-\int_{\bR^3}\p_z\f{(h^\th)^2}{r}w^\th dx-\int_{\bR^3}\p_r\rho w^\th dx\\
&:=I_1+I_2+I_3.
\eali
\ees
We estimate $I_i\ (i=1,2,3)$ separately. H\"{o}lder inequality and Gagliardo-Nirenberg interpolation inequality imply that
\bes
\bali
 I_1&\leq \pa u^r \pa_{L^3}\big\pa \f{w^\th}{r} \big\pa_{L^2}\pa w^\th \pa_{L^6}\\
    &\leq\pa u^r \pa_{L^3}\big\pa \O \big\pa_{L^2}\pa \na w^\th \pa_{L^2}\\
    &\leq C\pa u^r \pa^{2}_{L^3}\big\pa \O \big\pa^2_{L^2}+\f{1}{4}\pa \na w^\th \pa^2_{L^2}\\
    &\leq C\pa u^r \pa_{L^2}\pa\na u^r \pa_{L^2}\big\pa \O \big\pa^2_{L^2}+\f{1}{4}\pa \na w^\th \pa^2_{L^2},
\eali
\ees
and
\bes
\bali
 I_2&= \int_{\bR^3}\f{(h^\th)^2}{r}\p_zw^\th dx\\
    &\leq \pa  H\pa_{L^\i}\pa h^\th\pa_{L^2}\pa \na w^\th\pa_{L^2}\\
    &\leq C\pa  H\pa^2_{L^\i}\pa h^\th\pa^2_{L^2}+\f{1}{4}\pa \na w^\th\pa^2_{L^2}.
\eali
\ees
Also
\bes
\bali
 I_3&= -2\pi\int_{\bR}\int^\i_{0}\p_r\rho w^\th rdrdz\\
    &=2\pi\int_{\bR}\int^\i_{0}\rho \p_r(w^\th r)drdz\\
    &=2\pi\int_{\bR}\int^\i_{0}\rho \p_rw^\th rdrdz+\int_{\bR^3}\rho \f{w^\th}{r} dx\\
    &\leq \pa\rho\pa_{L^2}\pa\na w^\th\pa_{L^2}+\pa\rho\pa_{L^2}\big\pa \f{w^\th}{r}\big\pa_{L^2}\\
    &\leq C\pa\rho\pa^2_{L^2}+\f{1}{4}\Big(\pa\na w^\th\pa^2_{L^2}+\big\pa \f{w^\th}{r}\big\pa^2_{L^2}\Big).
\eali
\ees
The above estimates and Proposition \ref{pl2}, Proposition \ref{pomega} indicate that
\bes
\bali
&\q\f{d}{dt}\pa w^\th \pa^2_{L^2}+\pa \na w^\th \pa^2_{L^2}+\big\pa \f{w^\th}{r} \big\pa^2_{L^2}\\
&\leq C\pa u^r \pa_{L^2}\pa\na u^r \pa_{L^2}\big\pa \O \big\pa^2_{L^2}+C\pa  H\pa^2_{L^\i}\pa h\pa^2_{L^2}+C\pa\rho\pa^2_{L^2}\\
&\leq C_0(1+t)\Phi_{1,c_0}(t)\pa\na u^r \pa_{L^2}+C_0\pa H_0\pa^2_{L^\i}(1+t)^2+C\pa\rho_0\pa^2_{L^2}.
\eali
\ees
Integration on time implies that
\be\label{omegath}
\bali
&\q\pa w^\th(t) \pa^2_{L^2}+\int^t_0\pa \na w^\th(s) \pa^2_{ L^2}ds+\int^t_0\big\pa \f{w^\th}{r}(s) \big\pa^2_{L^2}ds\\
&\leq \Phi_{1,c_0}(t).
\eali
\ee
Using the identity $\na\times\na\times u=-\Dl u+\na\na\cdot u$ and divergence-free condition of $u$, we have
\be\label{bs}
\na u=\na (-\Dl)^{-1}\na\times w=\na (-\Dl)^{-1}\na\times(w^\th e_\th).
\ee
Calder\'{o}n-Zygmund theorem implies that for any $1<p<+\i$, we have
\be\label{cz}
\|\na u(t)\|_{L^p}\leq C_p\|w^\th(t)\|_{L^p},\qq \|\na^2 u(t)\|_{L^p}\leq C_p\lt(\|\na w^\th(t)\|_{L^p}+\lt\|\f{w^\th(t)}{r}\rt\|_{L^p}\rt).
\ee
From \eqref{omegath} and \eqref{cz}, we see that
\be\label{uh1}
\bali
&\q\pa \na u(t) \pa^2_{L^2}+\int^t_0\pa \na^2 u(s) \pa^2_{ L^2}ds\leq \Phi_{1,c_0}(t).
\eali
\ee

In order to bootstrap our energy estimates, we need the $L^1_tL^\i$ estimate of $u$. Before getting that, we first perform the $L^\i_tL^4$ estimates of $h^\th$ and $w^\th$.
\subsection{${L^\i_tL^4}$ estimate of ${h^\th}$ and ${w^\th}$ }

Performing $L^4$ inner product of $h^\th$ and using H\"{o}lder inequality, Gagliardo-Nirenberg interpolation inequality, we see that
\bes
\bali
\f{d}{dt}\pa h^\th(t)\pa^4_{L^4}&\leq 4\int_{\bR^3} \f{u^r}{r}(h^\th)^4dx\\
                                 &\leq 4\pa  H\pa_{L^\i}\int_{\bR^3} |u^r|(h^\th)^3 dx\\
                                 &\leq 4\pa  H_0\pa_{L^\i}\pa u^r\pa_{L^4}\pa h^\th\pa^3_{L^4}\\
                                 &\leq C\pa  H_0\pa_{L^\i}\pa \na u^r\pa^{3/4}_{L^2}\pa u^r\pa^{1/4}_{L^2}\pa h^\th\pa^3_{L^4}.
\eali
\ees
Integration on time implies that
\be\label{hl4}
\bali
\pa h^\th(t)\pa_{L^4}\leq \Phi_{1,c_0}(t).
\eali
\ee
Next performing the standard $L^4$ inner product of the $w^\th$ equation, we have
\bes
\bali
&\q\f{1}{4}\f{d}{dt}\pa w^\th \pa^4_{L^4}+\f{3}{4}\pa \na |w^\th|^2 \pa^2_{L^2}+\big\pa \f{|w^\th|^2}{r} \big\pa^2_{L^2}\\
&\leq \int_{\bR^3}\f{u^r}{r}(w^\th)^4dx-\int_{\bR^3}\p_z\f{(h^\th)^2}{r}(w^\th)^3 dx-\int_{\bR^3}\p_r\rho (w^\th)^3 dx\\
&:=I_1+I_2+I_3.
\eali
\ees
By the H\"{o}lder inequality, Gagliardo-Nirenberg interpolation inequality and Young inequality, we have
\bes
\bali
 I_1&\leq \pa u^r \pa_{L^4}\big\pa \f{w^\th}{r} \big\pa_{L^2}\pa (w^\th)^3 \pa_{L^4}\\
     &\leq C\pa u^r \pa^{1/4}_{L^2}\pa \na u^r \pa^{3/4}_{L^2}\big\pa \O \big\pa_{L^2}\pa (w^\th)^2 \pa^{3/2}_{L^6}\\
    &\leq C\pa u^r \pa^{1/4}_{L^2}\pa \na u^r \pa^{3/4}_{L^2}\big\pa \O \big\pa_{L^2}\pa \na(w^\th)^2 \pa^{3/2}_{L^2}\\
    &\leq C\pa u^r \pa_{L^2}\pa \na u^r \pa^{3}_{L^2}\big\pa \O \big\pa^4_{L^2}+\f{1}{8}\pa \na(w^\th)^2 \pa^{2}_{L^2}.
\eali
\ees
Also, H\"{o}lder inequality and Young inequality imply
\bes
\bali
 I_2&= \int_{\bR^3}\f{(h^\th)^2}{r}\p_z(w^\th)^3 dx\\
    &=3 \int_{\bR^3}\f{(h^\th)^2}{r} (w^\th)^2\p_zw^\th dx\\
    &\leq C\pa  H\pa_{L^\i}\pa h^\th\pa_{L^4}\pa w^\th\p_z w^\th\pa_{L^2}\pa w^\th\pa_{L^4}\\
    &\leq C\pa  H_0\pa^4_{L^\i}\pa h^\th\pa^4_{L^4}+\f{1}{8}\pa\p_z (w^\th)^2\pa^2_{L^2}+\pa w^\th\pa^4_{L^4},
\eali
\ees
and the same, we have
\bes
\bali
 I_3&= -2\pi\int_{\bR}\int^\i_{0}\p_r\rho (w^\th)^3 rdrdz\\
    &=2\pi\int_{\bR}\int^\i_{0}\rho \p_r\big((w^\th)^3 r\big)drdz\\
    &=6\pi\int_{\bR}\int^\i_{0}\rho (w^\th)^2\p_rw^\th rdrdz+\int_{\bR^3}\rho \f{(w^\th)^3}{r} dx\\
    &\leq C\pa\rho\pa_{L^\i}\pa\na (w^\th)^2\pa_{L^2}\pa w^\th\pa_{L^2}+\pa\rho\pa_{L^\i}\big\pa \f{(w^\th)^2}{r}\big\pa_{L^2}\pa w^\th\pa_{L^2}\\
    &\leq C\pa\rho\pa^2_{L^\i}\pa w^\th\pa^2_{L^2}+\f{1}{4}\pa\na (w^\th)^2\pa^2_{L^2}+\f{1}{4}\big\pa \f{(w^\th)^2}{r}\big\pa^2_{L^2}.
\eali
\ees
Using \eqref{uh1}, \eqref{hl4} and Proposition \ref{pl2}, the above inequalities imply
\bes
\bali
&\q\f{d}{dt}\pa w^\th \pa^4_{L^4}+\pa \na |w^\th|^2 \pa^2_{L^2}+\big\pa \f{|w^\th|^2}{r} \big\pa^2_{L^2}\\
&\leq C\pa w^\th\pa^4_{L^4}+C\pa u^r \pa_{L^2}\pa \na u^r \pa^{3}_{L^2}\big\pa \O \big\pa^4_{L^2}+C\pa  H_0\pa^4_{L^\i}\pa h^\th\pa^4_{L^4}+C\pa\rho\pa^2_{L^\i}\pa w^\th\pa^2_{L^2}\\
&\leq  C\pa w^\th\pa^4_{L^4}+\Phi_{1,c_0}(t).
\eali
\ees
Gronwall inequality implies that
\bes
\bali
\pa w^\th(t) \pa^4_{L^4}+\int^t_0\pa \na |w^\th(s)|^2 \pa^2_{L^2}ds+\int^t_0\big\pa \f{(w^\th)^2}{r}(s) \big\pa^2_{L^2}ds\leq \Phi_{1,c_0}(t).
\eali
\ees
The above inequality implies that
\be\label{uw14}
\pa \na u(t) \pa_{L^4}\leq \Phi_{1,c_0}(t).
\ee
Next we give a crucial estimate for bootstrapping the regularity of the solution.
\subsection{${L^1_tL^\i}$ estimate of ${\na u}$}

Applying $\na\times$ to $\eqref{MB1}_1$, we have
\be\label{we}
\p_t w-\Dl w=-\na\times[u\cdot\na u-h\cdot\na h-\rho e_3].
\ee
For a $H^1$ vector function $f$, we have
\bes
(\na\times f)\times f=f\cdot\na f-\f{1}{2}\na |f|^2.
\ees
Then we have
\bes
\na\times(f\cdot\na f)=\na\times[(\na\times f)\times f].
\ees
Inserting this into \eqref{we}, we have
\bes
\p_t w-\Dl w=-\na\times[(\na\times u)\times u-(\na\times h)\times h-\rho e_3].
\ees
Then we can write it as
\bes
\bali
w&=e^{t\Dl}w_0-\int^t_0e^{(t-s)\Dl}(\na\times[(\na\times u)\times u-(\na\times h)\times h-\rho e_3])ds\\
 &=e^{t\Dl}w_0-\int^t_0e^{(t-s)\Dl}\na\times[(\na\times u)\times u]ds\\
 &\q+\int^t_0e^{(t-s)\Dl}\na\times[(\na\times h)\times h]ds+\int^t_0e^{(t-s)\Dl}\na\times[\rho e_3]ds.
\eali
\ees
By a direct computation, if $h=h^\th e_\th$, we can get
\bes
\na\times[(\na\times h)\times h]=-2\f{h^\th}{r}\p_zh^\th e_\th=-\p_z(Hh^\th e_\th).
\ees
Then we have
\bes
\bali
w&=e^{t\Dl}w_0-\int^t_0e^{(t-s)\Dl}\na\times[(\na\times u)\times u]ds\\
 &\q-\int^t_0e^{(t-s)\Dl}\p_z( Hh^\th e_\th)ds+\int^t_0e^{(t-s)\Dl}\na\times[\rho e_3]ds.
\eali
\ees
Then by using \eqref{uw14}, the $L^s_tL^q_x$ $(1<s,q<+\i)$ estimates for the parabolic equation of singular integral and potentials (see, for example, \cite{Lje:1967PSPM,vWw:1982JLMS}) give that
\bes
\bali
&\pa\na w\pa_{L^2([0,t],L^4(\bR^3))}\\
\ls& \pa \na w_0\pa_{L^4(\bR^3)}t^{1/2}+\pa (\na\times u)\times u\pa_{L^2([0,t],L^4(\bR^3))}\\
&+\pa Hh^\th \pa_{L^2([0,t],L^4(\bR^3))}+\pa \rho \pa_{L^2([0,t],L^4(\bR^3))}\\
\ls&\pa \na w^\th_0\pa_{L^4(\bR^3)}t^{1/2}+\pa u\pa_{L^\i([0,t],L^\i(\bR^3))}\pa \na\times u\pa_{L^2([0,t],L^4(\bR^3))}\\
&+\pa  H \pa_{L^\i([0,t],L^\i(\bR^3))}\pa h^\th \pa_{L^2([0,t],L^4(\bR^3))}+\pa \rho_0 \pa_{L^2([0,t],L^4(\bR^3))}\\
\ls&\pa \na w_0\pa_{L^4(\bR^3)}t^{1/2}+\pa u\pa^{1/7}_{L^\i([0,t],L^2(\bR^3))}\pa\na u\pa^{6/7}_{L^\i([0,t],L^4(\bR^3))}\pa \na u\pa_{L^2([0,t],L^4(\bR^3))}\\
&+\pa H\pa_{L^\i([0,t],L^\i(\bR^3))}\pa h^\th \pa_{L^2([0,t],L^4(\bR^3))}+\pa \rho_0 \pa_{L^2([0,t],L^4(\bR^3))}\\
\leq& \Phi_{1,c_0}(t).
\eali
\ees
This, combining with \eqref{bs}, implies
\bes
\pa\na^2 u\pa_{L^2([0,t],L^4(\bR^3))}\leq C\pa\na w\pa_{L^2([0,t],L^4(\bR^3))}\leq \Phi_{1,c_0}(t).
\ees
Then by using H\"{o}lder inequality and  Gagliardo-Nirenberg interpolation inequality, we have
\be\label{duli}
\bali
\pa\na u\pa_{L^1([0,t],L^\i(\bR^3))}\leq& \int^t_0\pa \na u(s)\pa^{1/4}_{L^4}\pa \na^2 u(s)\pa^{3/4}_{L^4}ds\\
                                    \leq&\pa \na u(s)\pa^{1/4}_{L^\i[0,t],L^4(\bR^3)}\Big(\int^t_0\pa \na^2 u(s)\pa^2_{L^4}ds\Big)^{3/8}(\int^t_0ds)^{5/8}\\
                                    \leq& \Phi_{1,c_0}(t).
\eali
\ee
\begin{remark}
In cylindrical coordinates, for the axially symmetric velocity $u$, a direct computation indicates that
\be\label{cylind}
|\na u|\thickapprox |\t{\na} (u^r,u^\th, u^z)|+ \lt |\lt(\f{u^r}{r}, \f{u^\th}{r}\rt)\rt |,
\ee
where $\t{\na}=(\p_r,\p_z)$. From \eqref{duli} and \eqref{cylind}, we can also have
\be\label{cylind1}
\lt\pa\f{u^r}{r}\rt\pa_{L^1([0,t],L^\i(\bR^3))}\leq \Phi_{1,c_0}(t).
\ee
\end{remark}
\qed

Next we will use ${L^1_tL^\i}$ estimate of ${\na u}$ to bootstrap the regularity of the solution.
\subsection{ Estimate of ${\na \rho}$ and ${\na h}$}

Applying $\na$ to the third equation of \eqref{MB1}, we have
\bes
\p_t\na \rho+u\cdot\na \na\rho=-\na u\cdot\na \rho.
\ees
We can have for $1\leq p\leq 6$,
\bes
\bali
\pa \na \rho(t)\pa_{L^p}&\leq \pa \na \rho_0\pa_{L^p}+C\int^t_0\pa\na u\pa_{L^\i}\pa \na \rho(s)\pa_{L^p}ds.
\eali
\ees
Using the estimate \eqref{duli}, Gronwall inequality indicates that
\be\label{rhoh1}
\pa \na \rho(t)\pa_{L^p}\leq \Phi_{2,c_0}(t).
\ee
For the estimate of $\na h$, first we write the second equation of \eqref{MB1} as
\bes
\p_t h+u\cdot\na h=\f{u^r}{r}h.
\ees
Applying $\na$ to the above equality, we have
\bes
\p_t\na h+u\cdot\na \na h=-\na u\cdot\na h+\f{u^r}{r}\na h+\na u^r He_\th+(\na\f{1}{r})u^r h.
\ees
Noting
\bes
(\na\f{1}{r})u^r h=-\f{1}{r^2}e_ru^r h=-\f{u^r}{r} He_r\otimes e_\th,
\ees
and, as \eqref{cylind}, $| H|=|\f{h^\th}{r}|\ls |\na h| $,
we have, for $1\leq p\leq 6$,
\bes
\bali
\pa \na h(t)\pa_{L^p}&\leq \pa \na h_0\pa_{L^p}+C\int^t_0\pa(\na u,u^r/r)\pa_{L^\i}\pa \na h(s)\pa_{L^p}ds\\
&\q +C\int^t_0\pa(\na u,u^r/r)\pa_{L^\i}\pa H(s)\pa_{L^p}ds.
\eali
\ees
Also using the estimates \eqref{duli} and \eqref{cylind1}, Gronwall inequality indicates that
\be\label{hh1}
\pa \na h(t)\pa_{L^p}\leq \Phi_{2,c_0}(t).
\ee
Combining the estimates in \eqref{uh1}, \eqref{rhoh1} and \eqref{hh1}, we finish the proof of Proposition \ref{pH1} and \eqref{h1} is valid.
\section {${H^2}$ estimate of the solution and proof of Theorem \ref{th1}}

In this section, we give a prior $H^2$ estimate for the solution of system \ref{MB4}. We have the following Proposition.
\begin{proposition}\label{pH2}
Suppose $(u,h,\rho)$ be the smooth solution of \eqref{MB1} with initial data $(u_0,h_0,\rho_0)$ satisfying assumptions in Theorem \ref{th1}, then we have, for $t\in\bR_+$,
\be\label{h2}
\bali
&\q\pa (\na^2 u(t),\na^2 h(t),\na^2\rho(t)) \pa^2_{L^2}+\int^t_0\pa \na^3 u(s) \pa^2_{ L^2}ds\leq \Phi_{3,c_0}(t),
\eali
\ee
where $c_0$ is a positive constant depending only on $H^2$ norms of $u_0, h_0, \rho_0$ and $L^\i$ norm of $ H_0$.
\end{proposition}

\subsection{ Estimate of ${\na^2 u,\na^2 h}$}
Applying $\na^2$ to \eqref{MB1}, we have
\begin{equation}\label{MBSD}
\left\{
\begin{aligned}
\p_t\na^2u+u\cdot\nabla \na^2u+\nabla \na^2p-\Delta \na^2u-h\cdot\nabla \na^2h&=-[\na^2,u\cdot\na]u\\
&\qq+[\na^2,h\cdot\na]h+\na^2(\rho e_3),\\
\p_t \na^2h+u\cdot\nabla \na^2h-h\cdot\nabla \na^2u&=-[\na^2,u\cdot\na]h+[\na^2,h\cdot\na]u.\\
\end{aligned}
\right.
\end{equation}
Next we will use the following commutator estimate due to Kato-Ponce \cite{KP:1988CPAM},
\be\label{commutator}
\|\La^m(fg)-f\La^m g\|_{L^p}\leq C\Big(\|\na f\|_{L^{p_1}}\|\La^{m-1}g\|_{L^{p'_1}}+\|\La^m f\|_{L^{p_2}}\| g\|_{L^{p'_2}}\Big)
\ee
with $m\in\bN$, $\La=(-\Dl)^{1/2}$ and $1/p=1/p_1+1/p'_1=1/p_2+1/p'_2$.

Performing the $L^2$ energy estimate of \eqref{MBSD}, we have
\bes
\bali
&\q\f{1}{2}\f{d}{dt}\Big(\pa\na^2 u(t)\pa^2_{L^2}+\pa \na^2 h \pa^2_{L^2}\Big)+\pa\na^3u(t) \pa^2_{L^2}\\
&=- \int_{\bR^3}[\na^2,u\cdot\na]u \na^2 udx+\int_{\bR^3}[\na^2,h\cdot\na]h \na^2 udx-\int_{\bR^3}[\na^2,u\cdot\na]h \na^2 hdx\\
&\q +\int_{\bR^3}[\na^2,h\cdot\na]u \na^2 hdx+\int_{\bR^3}\na^2(\rho e_3)\na^2 udx\\
&:=I_1+I_2+I_3+I_4+I_5.
\eali
\ees
We estimate $I_i\ (i=1,2,3,4,5)$ term by term. Using \eqref{commutator}, Gagliardo-Nirenberg interpolation inequality and Young inequality, we have
\bes
\bali
I_1&\leq \pa[\na^2,u\cdot\na]u\pa_{L^2(\bR^3)}\pa\na^2u \pa_{L^2(\bR^3)}\\
&\leq\pa\na u\pa_{L^\i}\pa\na^2 u\pa_{L^2}\pa\na^2u \pa_{L^2(\bR^3)}\\
&\leq \pa\na u\pa_{L^\i}\pa\na^2 u\pa^2_{L^2},
\eali
\ees
and
\bes
\bali
I_2&\leq \pa[\na^2,h\cdot\na]h\pa_{L^{3/2}(\bR^3)}\pa\na^2u \pa_{L^3(\bR^3)}\\
&\leq\pa\na^2 u\pa_{L^3}\pa\na h\pa_{L^6}\pa\na^2h \pa_{L^2(\bR^3)}\\
&\leq\pa\na^2 u\pa^{1/2}_{L^2}\pa\na^3 u\pa^{1/2}_{L^2}\pa\na h\pa_{L^6}\pa\na^2h \pa_{L^2(\bR^3)}\\
&\leq C\pa\na h\pa^{4/3}_{L^6}(\pa\na^2 u\pa_{L^2}+\pa\na^2h \pa_{L^2(\bR^3)})^2+\f{1}{4}\pa\na^3 u\pa^2_{L^2}\\
&\leq\Phi_{2,c_0}(t)(\pa\na^2 u\pa_{L^2}+\pa\na^2h \pa_{L^2(\bR^3)})^2+\f{1}{4}\pa\na^3 u\pa^2_{L^2}.
\eali
\ees
Also the commutator estimate \eqref{commutator} implies
\bes
\bali
I_3&\leq \pa[\na^2,u\cdot\na]h\pa_{L^2(\bR^3)}\pa\na^2h\pa_{L^2(\bR^3)}\\
&\leq\big(\pa\na u\pa_{L^\i}\pa\na^2h \pa_{L^2(\bR^3)}+\pa\na^2 u\pa_{L^3}\pa\na h\pa_{L^6}\big)\pa\na^2h \pa_{L^2(\bR^3)}\\
&\leq\pa\na u\pa_{L^\i}\pa\na^2h \pa^2_{L^2(\bR^3)}+\Phi_{2,c_0}(t)(\pa\na^2 u\pa_{L^2}+\pa\na^2h \pa_{L^2(\bR^3)})^2+\f{1}{4}\pa\na^3 u\pa^2_{L^2}.
\eali
\ees
The same, we can get
\bes
\bali
I_4\leq\pa\na u\pa_{L^\i}\pa\na^2h \pa^2_{L^2(\bR^3)}+\Phi_{2,c_0}(t)(\pa\na^2 u\pa_{L^2}+\pa\na^2h \pa_{L^2(\bR^3)})^2+\f{1}{4}\pa\na^3 u\pa^2_{L^2},
\eali
\ees
and
\bes
\bali
|I_5|&\leq\Big|\int_{\bR^3}\na(\rho e_3)\na^3 udx\Big|\\
&\leq \pa\na \rho\pa_{L^2}\pa\na^3 u\pa_{L^2}\\
&\leq \f{1}{8}\pa\na^3 u\pa^2_{L^2}+C\pa\na \rho\pa^2_{L^2}.
\eali
\ees
The above estimates indicate that
\bes
\bali
&\q\f{1}{2}\f{d}{dt}\Big(\pa\na^2 u(t)\pa^2_{L^2}+\pa \na^2 h (t) \pa^2_{L^2}\Big)+\pa\na^3u(t) \pa^2_{L^2}\\
&\leq \big(\pa\na u (t)\pa_{L^\i}+\Phi_{2,c_0} (t)\big)\big(\pa\na^2 u (t)\pa^2_{L^2}+\pa\na^2 h (t)\pa^2_{L^2}\big)+\Phi_{2,c_0} (t).
\eali
\ees
Gronwall inequality indicates that
\be\label{uhh2}
\bali
\Big(\pa\na^2 u(t)\pa^2_{L^2}+\pa \na^2 h \pa^2_{L^2}\Big)+\int^t_0\pa\na^3u(s) \pa^2_{L^2}ds\leq \Phi_{3,c_0}(t).
\eali
\ee
\subsection{ Estimate of $\na^2 \rho$}

Next we give the estimate of $\na^2\rho$. Applying $\na^2$ to the third equation of \eqref{MB1}, we have
\bes
\p_t\na^2 \rho+u\cdot\na \na^2\rho=-[\na^2,u\cdot\na]\rho.
\ees
Standard $L^2$ energy estimate implies that
\bes
\bali
&\q\pa \na^2 \rho(t)\pa_{L^2}\\
&\leq \pa \na^2 \rho_0\pa_{L^2}+C\int^t_0\pa[\na^2,u\cdot\na]\rho\pa_{L^2}ds\\
&\leq \pa \na^2 \rho_0\pa_{L^2}+C\int^t_0\Big(\pa \na u\pa_{L^\i}\pa \na^2 \rho\pa_{L^2}+\pa \na^2 u\pa_{L^3}\pa \na \rho\pa_{L^6}\Big)ds\\
&\leq \pa \na^2 \rho_0\pa_{L^2}+C\int^t_0\Big(\pa \na u\pa_{L^\i}\pa \na^2 \rho\pa_{L^2}+\pa \na^2 u\pa^{1/2}_{L^2}\pa \na^3 u\pa^{1/2}_{L^2}\pa \na \rho\pa_{L^6}\Big)ds\\
&\leq \pa \na^2 \rho_0\pa_{L^2}+C\int^t_0\pa \na u\pa_{L^\i}\pa \na^2 \rho\pa_{L^2}ds+\Phi_{3,c_0}(t).
\eali
\ees
Gronwall inequality indicates that
\be\label{rhoh2}
\pa \na^2 \rho(t)\pa_{L^p}\leq \Phi_{3,c_0}(t).
\ee
The combination of \eqref{uhh2} and \eqref{rhoh2} proves Proposition \ref{pH2} and \eqref{h2} is valid.

\textbf{Proof of Theorem \ref{th1}.} Combining Proposition \ref{pl2}, Proposition \ref{pH1} and Proposition \ref{pH2}, we can get the a priori estimate \eqref{energy}. Then the local existence and uniqueness theorem in \cite{LP:JDE2017} and the a priori estimate \eqref{energy} together prove Theorem \ref{th1}.

\section*{Acknowledgments}

The author wish to thank Dr. Zijin Li in Nanjing University of Information Science and Technology for helpful discussions on this topic.

\qed

\bibliographystyle{plain}



\def\cprime{$'$}

\end{document}